\documentclass[12pt]{article}
\font\MSBMx			msbm10
\font\EUFMx			eufm10
   \font\EUFMvii			eufm7
    \font\EUFMv			eufm5
\usepackage{amssymb,latexsym,
amsfonts,amsbsy}

\newfam\German
\textfont\German=\EUFMx
\scriptfont\German=\EUFMvii
\scriptscriptfont\German=\EUFMv

\newfam\Blackboadbold
\textfont\Blackboadbold=\MSBMx
\scriptfont\Blackboadbold=\MSBMx

\newcommand{\qed}{\hbox{\rule[-2pt]
{3pt}{6pt}}}

\newtheorem{dfe}{Definition}
[section]
\newtheorem{lem}[dfe]{Lemma}
\newtheorem{theo}[dfe]{Theorem}

\newtheorem{pro}[dfe]{Proposition}

\makeatletter

\@addtoreset{equation}{section}
\makeatother



\begin{title}
{On antipodal spherical $t$-designs of degree $s$ with $t\geq 2s-3$}
\end{title}
\author{Eiichi Bannai and
Etsuko Bannai\\
\\
Graduate School of Mathematics
\\ Kyushu University}
\date{}
\begin{document}
\maketitle

\begin{abstract}
We prove that if $X$ is a spherical $t$-design and $s$-distance set with $t\geq 2s-3$, 
then $X$ has the structure of Q-polynomial association scheme of class $s$. Also, we describe the parameters of the association scheme. 
\end{abstract}

\section{Introduntion}
Delsarte-Goethals-Seidel \cite{D-G-S} studied finite subsets $X$ in the 
unit sphere $S^{n-1}=\{x=(x_1,x_2,\ldots ,x_n)\in \mathbb R^n\mid 
x_1^2+x_2^2+\cdots +x_n^2=1 \}$ 
from the viewpoint of algebraic 
combinatorics. Here the following two parameters $s=s(X)$ and $t=t(X)$ 
play important roles.  The $s$ is called the degree of $X$ and 
is defined as the number of distinct distances between two 
distinct elements of $X$, (then $X$ is called an $s$-distance set), 
while  $t$ is called the strength of $X$ and is usually defined as the 
largest $t$ such that $X$ becomes a $t$-design. Note that when we 
say $X$ is a $t$-design, the strength of $X$ might be larger than 
$s(X)$. So our use of $t$ is slightly ambiguous, but this convention 
is very useful and we believe that serious confusion 
will not occur. The important results due to 
Delsarte-Goethals-Seidel
\cite{D-G-S} are as follows 
(the leader is referred to \cite{D-G-S,B-B-1} for the 
definition of undefined terminologies): \\
(i) We always have $t\leq 2s.$\\
(ii) If $t\geq s-1,$ then $X$ is distance invariant. \\
(iii) If $t\geq 2s-2,$ then $X$ has the structure of 
Q-polynomial 
association scheme of class $s$. \\
(iv) $t=2s$ if and only if $X$ is a tight $2s$-design. \\
(v) $t=2s-1$ and $X$ is antipodal, i.e., if $x\in X$ then 
$-x\in X,$ if and only if $X$ is a tight $(2s-1)$-design. \\

Note that tight $t$-designs are classified except for 
$t\not= 4,5, 7$ (see \cite{B-D-1,B-D-2,B-S})
and there are some 
recent development on the study of tight $t$-designs 
with $t=4,5,7$ (cf. \cite{B-M-V}.) So, it is important to study, or to classify, 
(or to show the nonexistence for the case of large $t$ and $n\geq 3$), $X$ with $t$ close to $2s$. \\
In this paper, we study the case $t=2s-3$ with the additional 
assumption that $X$ is antipodal, and show that similar 
properties for the general non-antipodal case with $t\geq 2s-2$ 
hold true (cf. (iii) above). Namely, in this paper, we prove the following theorems. 

\begin{theo} Let $X$ be a spherical $t$-design
and $s$-distance set. Assume $X$ is antipodal
and $t\geq 2s-3$. Then $X$ has the structure 
of an association scheme of class $s$. 
\label{theo:main1}
\end{theo}
\begin{theo}\label{theo:main2} Let $P=(P_j(i))$, $Q=(Q_j(i))$ be the
first and second eigen matrices
and $q_{i,j}^k,\ 0\leq i,j,k\leq s$, be the
Krein numbers of the association scheme
given in Theorem 1.1. 
  Then by a suitable ordering of the 
  adjacency matrices $D_0,\ D_1,\ \ldots,\ D_s$
  and the
  basis of primitive idempotents $E_0,E_1,\ldots, E_s$
  we obtain the followings.
\begin{enumerate}
\item $P_1(i)=(-1)^i$ for any $i$ with $0\leq i\leq s$.
\item $Q_j(2i+1)=(-1)^jQ_j(2i)$ for any $i$ and $j$ wuth $1\leq 2i+1\leq s$ and $0\leq j\leq s$.
If $s$ is even, then $Q_j(s)=0$ for any odd
integer $j$ satisfying $1\leq j\leq s-1$.
\item $m_i=Q_i(0)={n+i-1\choose i}-{n+i-3\choose s-2}$
for $0\leq i\leq s-2$,

 $m_{s-1}=Q_{s-1}(0)=
\frac{|X|}{2}-{n+s-4\choose s-3}
$ 
 and
 
 $m_s=\frac{|X|}{2}-{n+s-3\choose s-2}$.
 
\item $q_{\mu,i}^j=0$ for any $\mu,i,j$
satisfying $0\leq  \mu,i,j\leq s$ and
$\mu+i+j=0$.
\item The dual intersection matrix $B_1^*=
(q_{1,i}^j)$ is tri-diagonal, that is 
the association scheme is a
Q-polynomial scheme.
More precisely 
$$B_1^*=\left[\begin{array}{ccccc}
*&c^*_1&\cdots&c_{s-1}^*&n\\
0&0&\cdots&\cdots&0\\
n&b_1^*&\cdots&b_{s-1}^*&*
\end{array}\right]$$
where 

$c_i^*=q_{1,i-1}^i=\frac{ni}{n+2i-2}$, 
for $i=1,\ldots,s-2$,

$b_i^*=q_{1,i}^{i-1}=\frac{n(n+i-3)}{n+2i-4}$ for $i=1,\ldots,s-3$,

$c_{s-1}^*=q_{1,s-2}^{s-1}=
 \frac{2n(n-1)(n+s-4)!}{(s-2)!(n-1)!|X|-2(s-2)(n+s-4)!}
$

$b_{s-2}^*=q_{1,s-1}^{s-2}=\frac{n(n+s-4)}
{n+2s-6}$,
 
 and
  
$b_{s-1}^*=q_{1,s}^{s-1}=
 \frac{(s-2)!n!|X|-2n(n+s-3)!}{(s-2)!(n-1)!|X|-2(s-2)(n+s-4)!} 
$.
\end{enumerate}
\end{theo}

In \S 2, we give some basic facts on
spherical $t$-designs and $s$-distance sets.
In \S3, we prove that $X$ has the structure 
of an association scheme of class $s$.
In \S4 we prove Theorem 1.2. 

\section{Some basic facts about spherical
$t$-designs}
We use the method given in the
paper by Delsarte-Goethals-Seidel \cite{D-G-S} to 
prove our main results. First we introduce
some definition and notations.
We denote by $\tilde Q_l(x)$ the Gegenbauer  
polynomial of degree $l$ attached to 
the unit sphere $S^{n-1}\subset \mathbb R^n$.
We use the notation $\tilde Q$ in order to 
distinguish it from the second eigenmatrix
$Q$.   
Here we use the normalization so that
$\tilde Q_l(1)=\dim(\mbox{Harm}_l(\mathbb R^n))$
holds. It is well known that the Gegenbauer polynomials 
satisfy the following equations (see \cite{D-G-S, B-B-1}).
\begin{eqnarray}
 &&\frac{l+1}{n+2l}\tilde Q_{l+1}(x)=x\tilde Q_l(x)-
 \frac{n+l-3}{n+2l-4}\tilde Q_{l-1}(x)
\quad  \mbox{for any integer $l\geq 0$ and}\\
 &&\tilde Q_i(x)\tilde Q_j(x)
=\sum_{{k=|i-j|}\atop
k\equiv i+j \ (\mbox{\scriptsize mod}\ 2)}^{i+j}
q_k(i,j)\tilde Q_k(x)\quad \mbox{for any 
integer $i,j,k\geq 0$},
 \end{eqnarray} 
where $q_k(i,j)$ is a nonnegative
real number. 
Let 
\begin{equation}x^\lambda=\sum_{l=0}^\lambda f_{\lambda,l}
\tilde Q_l(x)\label{equ:xlamda}\end{equation}
be the
Gegenbauer expansion of $x^\lambda$
for any non negative integer $\lambda$.
For each pair $\lambda,\mu$
of nonnegative integers, we define a polynomial by
$F_{\lambda,\mu}(x)=\sum_{l=0}^{\min\{\lambda,\mu\}}
f_{\lambda,l}f_{\mu,l}\tilde Q_l(x)$.
We denote by $\boldsymbol x\cdot \boldsymbol y$
the canonical inner product
between the
vectors $\boldsymbol x, \boldsymbol y\in
\mathbb R^n$.
For a finite subset $X\subset S^{n-1}$ we define
$A(X)=\{\boldsymbol x\cdot \boldsymbol y\ |\
\boldsymbol x, \boldsymbol y\in X
,\boldsymbol x\not= \boldsymbol y\}$
and $A'(X)=A(X)\cup \{1\}$.
Let $s=|A(X)|$. Then $X$ is called an
$s$-distance set.
Now we express
$A(X)=\{\alpha_1,\ 
\alpha_2,\ldots,\alpha_s\}$. 
Let
$\alpha_0=1$. 
We define matrices $D_i$ indexed by $X$
by 
\begin{equation}
D_i(\boldsymbol x, \boldsymbol y)=\left\{
\begin{array}{ll}1&\mbox{if $\boldsymbol x\cdot\boldsymbol y=\alpha_i$},\\
0&\mbox{if $\boldsymbol x\cdot\boldsymbol y
\not=\alpha_i$}.
\end{array}\right.
\end{equation}
Let $\mbox{Harm}(\mathbb R^n)$ be the
vector space of the harmonic polynomials
on $n$ variables 
$\boldsymbol x=(x_1,x_2,\ldots,x_n)$. 
We define a positive definite innerproduct
on $\mbox{Harm}(\mathbb R^n)$ 
by
$$\langle\varphi,\psi\rangle=\frac{1}{|S^{n-1}|}\int_{S^{n-1}}
\varphi(\boldsymbol x)\psi(\boldsymbol x)
d\sigma(\boldsymbol x)
$$
for $\varphi,\psi\in\mbox{Harm}(\mathbb R^n)$. 
Let $\mbox{Harm}_l(\mathbb R^n)$ 
be the subspace of 
$\mbox{Harm}(\mathbb R^n)$ generated by
homogeneous harmonic polynomial
 of degree $l$.
 Let $\varphi_{l,1},
 \ldots,\varphi_{l,h_l}$
be an orthonormal basis of $\mbox{Harm}_l(\mathbb R^n)$ with respect to the inner product given above.
The following addition
formula is well known.

\noindent
{\bf Addition formula}
\begin{eqnarray}\label{equ:addition}
\sum_{i=1}^{h_l}\varphi_{l,i}(\boldsymbol x)
\varphi_{l,i}(\boldsymbol y)
=\tilde Q_l(\boldsymbol x\cdot \boldsymbol y)
\end{eqnarray} 
holds for any 
$\boldsymbol x, \boldsymbol y\in
S^{n-1}$.\\

Let $H_l$ be the matrix
indexed by $X$ and 
$\{\varphi_{i,l}\mid 1\leq i\leq h_l\}$,
whose $(x,\varphi_{l,i})$-entry
is $\varphi_{l,i}(\boldsymbol x)$.
The following proposition is well known.
\begin{pro}
\label{pro:idmpotents} Let $X$ be a spherical $t$-design
and $s$-distance set. The notation is 
given as before. 
\begin{enumerate}
\item ${^tH_k}H_l=\delta_{k,l}I$ for any
nonnegative integers $k,\ l$ satisfying $0\leq k+l\leq t$.
\item $H_k\ {^tH_k}=\sum_{i=0}^s
\tilde Q_k(\alpha_i)D_i$ for any nonnegative integer
$k$.
\item $H_k\ {^tH_k}\neq 0$ for $k=0,1,\ldots s$.
\end{enumerate} 
\end{pro}
{\bf Proof} 
Let us compute the $(\varphi_{k,i},\varphi_{l,j})$-entry
of ${^tH_k}H_l$. Since $X$ is a spherical $t$-design,
we obtain
$$\sum_{\boldsymbol x\in X}\varphi_{k,i}
(\boldsymbol x)\varphi_{l,j}
(\boldsymbol x)=\frac{|X|}{|S^{n-1}|}
\int_{S^{n-1}}\varphi_{k,i}(\boldsymbol x)
\varphi_{l,j}(\boldsymbol x)d\sigma(\boldsymbol x)
=|X|\delta_{k,l}\delta_{i,j}.$$
This implies (1).
Let us compute the 
$(\boldsymbol x,\boldsymbol y)$-entry
of $H_k\ {^tH_k}$. Then addition formula 
(\ref{equ:addition}) 
implies
$$\sum_{i=1}^{h_k}\varphi_{i,k}(\boldsymbol x)
\varphi_{i,k}(\boldsymbol y)
=\tilde Q_k(\boldsymbol x\cdot\boldsymbol y).$$
This implies (2).
Since $\tilde Q_k(x)$ is a polynomial of degree
$k\leq s$ and $\alpha_0,\alpha_1,\ldots,\alpha_s$
are distinct $s+1$ real numbers,
we nust have an $\alpha_i$
satisfying $\tilde Q_k(\alpha_i)\neq 0$.
This implies (3). \hfill\qed\\

\section{Proof of Theorem 1.1}

\noindent
Let $X$ be a spherical 
$t$-design and  an $s$-distance set.
assume $X$ is antipodal.
Then, if $s$ is odd, then by arranging the
numbering of the elements in
$A(X)=\{\alpha_i\mid 1\leq i\leq s\}$
we may assume $\alpha_1=-1$,
$\alpha_{2i+1}=-\alpha_{2i}\neq 0$
for $i=0,\ \ldots, \frac{s-1}{2}$.
If $s$ is even,
we may assume $\alpha_1=-1$,
$\alpha_{2i+1}=-\alpha_{2i}\neq 0$
for $i=0,\ \ldots, \frac{s}{2}-1$
and $\alpha_s=0$. Note that
$\alpha_0=1$.
For $\boldsymbol x,\boldsymbol y\in X$
and $\alpha,\beta\in A'(X)$,
we define 
$$p_{\alpha,\beta}
(\boldsymbol x,\boldsymbol y)
=|\{\boldsymbol z\mid 
\boldsymbol x\cdot\boldsymbol z=\alpha\
\mbox{and}\ 
\boldsymbol z\cdot\boldsymbol y=\beta\}|.$$

\begin{lem}\label{lem:3-1} 
Definitions and notation are given as before.
Let $\boldsymbol x\cdot\boldsymbol y=\gamma$.
Then
 $p_{\alpha,\beta}
(\boldsymbol x,\boldsymbol y)$ depends
only on $\alpha,\beta$ and $\gamma$
and does not depend on the choice of 
$\boldsymbol x,\boldsymbol y$
satisfying $\boldsymbol x\cdot\boldsymbol y
=\gamma$.
\end{lem}
{\bf Proof} 
We have the following
\begin{eqnarray}&&p_{\alpha_i,\alpha_0}
(\boldsymbol x,\boldsymbol y)=
p_{\alpha_i,\alpha_0}(\boldsymbol x,\boldsymbol y)
=\left\{
\begin{array}{ll}1&\mbox{if $\alpha_i=\gamma$,}\\
0&\mbox{otherwise.}\\
\end{array}\right.\\
&&p_{\alpha_i,\alpha_1}
(\boldsymbol x,\boldsymbol y)=
p_{\alpha_i,\alpha_1}(\boldsymbol x,\boldsymbol y)
=\left\{
\begin{array}{ll}1&\mbox{if $\alpha_i=-\gamma$,}\\
0&\mbox{otherwise.}\\
\end{array}\right.
\end{eqnarray}
Let $\boldsymbol x,\boldsymbol y\in X$
and $\boldsymbol x\cdot\boldsymbol y=\gamma$.
Then for any $\lambda$ and $\mu$
satisfying $0\leq \lambda,\mu\leq s-2$ we compute
the $(\boldsymbol x,\boldsymbol y)$-entry
of the matrix
$$(\sum_{k=0}^\lambda f_{\lambda,k}H_k\ {^tH_k})
(\sum_{l=0}^\mu f_{\mu,l}H_l\ {^tH_l})$$
in two defferent ways,
where $ f_{\lambda,k}(0\leq k\leq \lambda),\ 
f_{\mu,l}(0\leq l\leq \mu)$ are defined in 
(\ref{equ:xlamda}).
First use Proposition \ref{pro:idmpotents} (1), 
and then use Proposition \ref{pro:idmpotents} 
(2), then we obtain
\begin{eqnarray}
&&((\sum_{k=0}^\lambda f_{\lambda,k}H_k\ {^tH_k})
(\sum_{l=0}^\mu f_{\mu,l}H_l\ {^tH_l}))
(\boldsymbol x,\boldsymbol y)
=|X|\sum_{k=0}^{\min\{\lambda,\mu\}} 
f_{\lambda,k}f_{\mu,k}(H_k\ {^tH_k})(\boldsymbol x,\boldsymbol y)\nonumber\\
&&=|X|\sum_{k=0}^{\min\{\lambda,\mu\}} 
f_{\lambda,k}f_{\mu,k}\tilde Q_k(\gamma)=|X|F_{\lambda,\mu}
(\gamma).
\label{equ:3-3}
\end{eqnarray}
Next, first apply Proposition \ref{pro:idmpotents} (2) and
then use (\ref{equ:xlamda}). Then we obtain
\begin{eqnarray}
&&((\sum_{k=0}^\lambda f_{\lambda,k}H_k\ {^tH_k})
(\sum_{l=0}^\mu f_{\mu,l}H_l\ {^tH_l}))(\boldsymbol x,\boldsymbol y)\nonumber\\
&&=((\sum_{k=0}^\lambda f_{\lambda,k}\sum_{i=0}^s
Q_k(\alpha_i)D_i)
(\sum_{l=0}^\mu f_{\mu,l}\sum_{j=0}^s
Q_l(\alpha_j)D_j))(\boldsymbol x,\boldsymbol y)
\nonumber\\
&&=\sum_{k=0}^\lambda\sum_{l=0}^\mu
\sum_{i=0}^s\sum_{j=0}^sf_{\lambda,k}f_{\mu,l}
Q_k(\alpha_i)Q_l(\alpha_j)
\sum_{\boldsymbol z\in X}
D_i(\boldsymbol x,\boldsymbol z)
D_j(\boldsymbol z,\boldsymbol y)
\nonumber\\
&&=\sum_{i=0}^s\sum_{j=0}^s
\sum_{k=0}^\lambda\sum_{l=0}^\mu
f_{\lambda,k}Q_k(\alpha_i)f_{\mu,l}Q_l(\alpha_j)
p_{\alpha_i,\alpha_j}(\boldsymbol x,\boldsymbol y)\nonumber\\
&&=\sum_{i=0}^s\sum_{j=0}^s\alpha_i^\lambda
\alpha_j^\mu p_{\alpha_i,\alpha_j}(\boldsymbol x,\boldsymbol y)=
\sum_{i=2}^s\sum_{j=2}^s\alpha_i^\lambda
\alpha_j^\mu p_{\alpha_i,\alpha_j}(\boldsymbol x,\boldsymbol y)
+\sum_{i=0}^1\sum_{j=2}^s
\alpha_i^\lambda
\alpha_j^\mu p_{\alpha_i,\alpha_j}(\boldsymbol x,\boldsymbol y)\nonumber\\
&&+\sum_{i=2}^s\sum_{j=0}^1
\alpha_i^\lambda
\alpha_j^\mu p_{\alpha_i,\alpha_j}(\boldsymbol x,\boldsymbol y)
+\sum_{i=0}^1\sum_{j=0}^1
\alpha_i^\lambda
\alpha_j^\mu p_{\alpha_i,\alpha_j}(\boldsymbol x,\boldsymbol y)\nonumber\\
&&=\sum_{i=2}^s\sum_{j=2}^s\alpha_i^\lambda
\alpha_j^\mu p_{\alpha_i,\alpha_j}(\boldsymbol x,\boldsymbol y)\nonumber\\
&&\qquad +\left\{
\begin{array}{ll}
\gamma^\mu +(-1)^\lambda(-\gamma)^\mu
 +\gamma^\lambda+(-\gamma)^\lambda(-1)^\mu &
 \mbox{if $\gamma\neq \alpha_0,\alpha_1$}\\
   \delta_{\gamma,\alpha_0}
+
(-1)^\mu \delta_{\gamma,\alpha_1}+
(-1)^\lambda \delta_{\gamma,\alpha_1}
+
(-1)^{\lambda+\mu}
\delta_{\gamma,\alpha_0}
&\mbox{if $\gamma\in\{\alpha_0,\alpha_1\}$}
 \end{array}
 \right.
 \label{equ:3-4}
\end{eqnarray}
Thus (\ref{equ:3-4}) and (\ref{equ:3-4}) yield a system of
linear equation whose indeterminates 
are $p_{\alpha_i,\alpha_j}(\boldsymbol x,\boldsymbol y)$,
$2\leq i,j\leq s$ and the coefficient matrix
is $W\otimes W$, where
$$W=\left[\begin{array}{cccc}
1&1&\cdots&1\\
\alpha_2&\alpha_3&\cdots&\alpha_s\\
\alpha_2^2&\alpha_3^2&\cdots&\alpha_s^2\\
\vdots&\vdots&\vdots&\vdots\\
\alpha_2^{s-2}&\alpha_3^{s-2}&\cdots&
\alpha_s^{s-2}
\end{array}
\right]$$
Since $W$ is invertible, 
$p_{\alpha_i,\alpha_j}
(\boldsymbol x,\boldsymbol y),\ (2\leq i,\ j \leq s)$
are uniquely determined by
$\alpha_i,\alpha_j$ and $\gamma$
and does not depend on the 
choice of $\boldsymbol x,\boldsymbol y$ satisfying 
$\boldsymbol x\cdot\boldsymbol y=\gamma$.
This completes the proof.
\hfill\qed\\

Lemma \ref{lem:3-1} implies that the spherical $t$-design
given in Theorem \ref{theo:main1} has the structure of an 
association scheme.\\

\section{Proof of Theorem 1.2}
First we will show the following proposition. 

\begin{pro}\label{pro:3-2}
The definition and notation are given as before.  
Let $F_i=\frac{1}{|X|}H_i\ {^tH}_i$ 
for $i=0,\ldots, s-2$.
Then $\{F_0,\ F_1,\ \ldots,\ F_{s-2}\}$ is a subset of the basis of the primitive idempotents
of the Bose-Mesner algebra
$\mathfrak A=\langle D_0,D_1,\ldots,D_s\rangle$.
\end{pro}
{\bf Proof} By definition $F_0=\frac{1}{|X|}J=E_0$,
where $J$ is the matrix whose elements are
all $1$. Proposition \ref{pro:idmpotents} implies 
that $F_i\in\mathfrak A$ and $F_iF_j=\delta_{i,j}F_i$
for any $0\leq i,j\leq s-2$.
Let 
$F_{s-1}=\frac{1}{|X|}H_{s-1}{^tH}_{s-1}$
and $F_s=I-\sum_{i=0}^{s-2}F_i$.
Then $F_{s-1},F_s\in \mathfrak A$.
Moreover $F_{s-1}F_i=F_sF_i=0$ and 
$F_s^2=F_s$
holds for any $i=0,\ \ldots,\ s-2$.
Since $F_0, F_1,\ldots,F_{s-2}$ are
idempotents, each $F_i$ is a partial some
of the primitive idempotents
$E_0,E_1,\ldots,E_s$ of $\mathfrak A$.
Assume that there exists an $F_i\not\in
\{E_0,E_1,\ldots,E_s\}$, then we have
decomposition $F_i=E+E'$, satisfying $E,E'\neq 0$,
$E^2=E,\ (E')^2=E'$ and $EE'=0$.
Then, $F_iE=E$ and $F_iE'=E'$ implies
$F_{s-1}E=F_{s-1}E'$ and $F_sE=F_sE'=0$.
On the other hand $F_0,F_1,\ldots,F_{i-1},E,E',F_{i+1},
\ldots, F_{s-2}$ are linearly independent.
Since Proposition \ref{pro:idmpotents} (3) implies $F_{s-1}\neq 0$, $\dim(\langle F_{s-1},F_s\rangle)=1$.
Hence there exists a real number $c$
satisfying $F_s=cF_{s-1}$.
This implies
$$
D_0-\frac{1}{|X|}\sum_{i=0}^{s-2}
\sum_{j=0}^{s}\tilde Q_{i}(\alpha_j)
D_j=c\frac{1}{|X|}\sum_{j=0}^s\tilde Q_{s-1}(\alpha_j)
D_j.$$
Then we obtain
$$c\tilde Q_{s-1}(\alpha_j)+\sum_{i=0}^{s-2}\tilde Q_i(\alpha_j)=0 $$
for any $j=1,\ \ldots\ ,\ s$.
However $c\tilde Q_{s-1}(x)+\sum_{i=0}^{s-2}
\tilde Q_i(x)$ is a polynomial of 
degree at most $s-1$ and $\alpha_1,\ldots,\alpha_s$ are
distinct to each other, this is a contradiction.
Hence $F_i\in \{E_j\mid 0\leq j\leq s\}$ for any 
$0\leq i\leq s-2$.  
\hfill\qed\\

\noindent
{\bf Proof of Theorem 1.2 (1)}\\ 
 Let $E_0,E_1,\ldots,E_d$
 be the basis of primitive idempotents
 of $\mathfrak A$ satisfying $E_i=\frac{1}{|X|}H_i\ {^tH}_i$
 for $i=0,1,\ldots,s-2$.
We note that we defined $\alpha_0=1$, $\alpha_1=-1$,
$\alpha_{2i+1}=-\alpha_{2i}$ for $0\leq 2i+1\leq s$,
and if $s$ is even $\alpha_s=0$.
Let $B_1$ be the intersection matrix of $\mathfrak A$ whose $(i,j)$-entry is
defined by
$B_1(i,j)=p_{1,i}^j$.
Then we have
\begin{eqnarray}
B_1=\left[
\begin{array}{cccc}
S&0&\cdots&0\\
0&S&\ddots&\vdots\\
\vdots&\ddots&\ddots&0\\
0&\cdots&0&S
\end{array}\right]\quad \mbox{if $s$ is odd,}
\quad
B_1=\left[
\begin{array}{cccc}
S&0&\cdots&0\\
0&\ddots&\ddots&\vdots\\
\vdots&\ddots&S&0\\
0&\cdots&0&1
\end{array}\right]\quad \mbox{if $s$ is even}.
\label{equ:B10}
\end{eqnarray}
where $S=\left[\begin{array}{cc}0&1\\1&0
\end{array}\right]$.
Hence the matrix $B_1$ has the
eigenvalue $1$ with multiplicity $[\frac{s+2}{2}]$ 
and $-1$ with multiplicity $[\frac{s+1}{2}]$.
Let us compute $D_1D_j$ in two different ways.
First we have
$$D_1D_j=\sum_{l=0}^sp_{1,j}^lD_l
=\sum_{l=0}^sp_{1,j}^l
\sum_{\mu=0}^sP_l(\mu)E_\mu
=\sum_{\mu=0}^s(\sum_{l=0}^sp_{1,j}^l
P_l(\mu))E_\mu.$$
On the other hand
$$D_1D_j=(\sum_{\mu=0}^sP_1(\mu)E_\mu)
(\sum_{\lambda=0}^sP_j(\lambda)E_\lambda)
=\sum_{\mu=0}^sP_1(\mu)P_j(\mu)E_\mu
$$
Hence we obtain
\begin{equation}\sum_{l=0}^s
p_{1,j}^lP_l(\mu)
=P_j(\mu)P_1(\mu)
\label{equ:B1}
\end{equation}
for any $j,\ \mu=0,1,\ldots, s$.
Let $M_1$ be a diagonal matrix 
whose diagonal entries are defined by
$M_1(i,i)=P_1(i)$ for $i=0,1,\ldots,s$.  
Then (\ref{equ:B1}) implies
\begin{equation}
B_1\ {^tP}={^tP} M_1.
\end{equation}
Hence
${^tP}M_1(^tP)^{-1}= B_1$ holds. This implies
that $\{P_1(0),P_1(1),\ldots, P_1(s)\}$
is the set of eigenvalues of $B_1$.
Hence $P_1(i)= 1$ or $-1$ for
any $0\leq i\leq s$ and 
$|\{i\mid  P_1(i)= 1\}|=[\frac{s+2}{2}]$
and $|\{i\mid  P_1(i)= -1\}|=[\frac{s+1}{2}]$
holds. 
On the other hand, 
$E_i=\frac{1}{|X|}H_i\ {^tH}_i
=\frac{1}{|X|}\sum_{j=0}\tilde Q_i(\alpha_j)D_j$ for $0\leq i\leq s-2$, and (\ref{equ:B10}) imply
\begin{eqnarray}
D_1E_i&=&\frac{1}{|X|}\sum_{j=0}^s
\tilde Q_i(\alpha_j)
D_1D_j
=\frac{1}{|X|}\sum_{l=0}^s
\left(\sum_{j=0}^s \tilde Q_i(j)p_{1,j}^l\right)D_l\nonumber\\
&=&\frac{1}{|X|}\sum_{l=0}^s
\left(\sum_{j=0}^s\tilde Q_i(j)B_1(l,j)\right)D_l
=\frac{1}{|X|}\sum_{l=0}^s(B_1\tilde Q)(l,i)D_l
\nonumber\\
&=&\frac{1}{|X|}\left(\tilde Q_i(\alpha_1)D_0
+\tilde Q_i(\alpha_0)D_1+\ldots+
\tilde Q_i(\alpha_{2j+1})D_{2j}
+\tilde Q_i(\alpha_{2j})D_{2j+1}+\cdots
\right)\nonumber\\
&=&\frac{1}{|X|}\sum_{j=0}^s\tilde Q_i(-\alpha_j)D_j
\nonumber\\
&=&(-1)^iE_i,\\
D_1E_i&=&\sum_{l=0}^sP_1(l)E_lE_i
=P_1(i)E_i.
\end{eqnarray}
Hence we obtain
$P_1(i)=(-1)^i$ for $0\leq i\leq s-2$.
Hence by arranging the ordering
of $E_{s-1}$ and $E_s$, we obtain
$P_1(i)=(-1)^i$ for any $0\leq i\leq s$.
This completes the proof.\\

\noindent
{\bf Proof of Theorem 1.2 (2)} \\
If $0\leq j\leq s-2$, then 
$Q_j(l)=\tilde Q_j(\alpha_l)$.
 By definition $\alpha_{2i+1}=-\alpha_{2i}$
for $0\leq i\leq \frac{s-1}{2}$.
Hence $Q_j(2i+1)=\tilde Q_j(\alpha_{2i+1})
=\tilde Q_j(-\alpha_{2i})=(-1)^j\tilde Q_j(\alpha_{2i})
=(-1)^j Q_j(2i)$.
If $s$ is even, then $\alpha_s=0$.
Hence $Q_j(s)=\tilde Q_j(\alpha_s)=\tilde Q_j(0)=0$
for any odd integer $j$ satisfying
$0\leq j\leq s-2$.
 Let us consider the condition
$QP=|X|I$. If $1\leq 2i+1\leq s$, we have
the following equations.
\begin{eqnarray}&&
\sum_{j=0}^sQ_j(2i)=\delta_{i,0}|X|,
\label{equ:ee}
\\
&&
\sum_{j=0}^sQ_j(2i)P_1(j)=\sum_{j=0}^s(-1)^jQ_j(2i)=0,
\label{equ:eo}\\
&&\sum_{j=0}^sQ_j(2i+1)P_0(j)
=\sum_{j=0}^{s-2}(-1)^jQ_j(2i)+Q_{s-1}(2i+1)
+Q_s(2i+1)=0,\label{equ:oe}\\
&&\sum_{j=0}^sQ_j(2i+1)P_1(j)
=\sum_{j=0}^s(-1)^jQ_j(2i+1)
\nonumber\\
&&=\sum_{j=0}^{s-2}Q_j(2i)+(-1)^{s-1}Q_{s-1}(2i+1)
+(-1)^sQ_{s}(2i+1)
=\delta_{i,0}|X|,\label{equ:oo}
\label{equ:3-14}
\end{eqnarray}
hold. Then (\ref{equ:eo}) and (\ref{equ:oe})
imply 
\begin{equation}
(-1)^{s-1}(Q_{s-1}(2i)-Q_s(2i))
=Q_{s-1}(2i+1)+Q_s(2i+1).
\label{equ:3-15}\end{equation}
Then  (\ref{equ:ee}) and (\ref{equ:oo})
imply 
\begin{equation}Q_{s-1}(2i)+Q_s(2i)
=(-1)^{s-1}(Q_{s-1}(2i+1)-Q_s(2i+1)).
\end{equation}
(\ref{equ:3-14}) and (\ref{equ:3-15}) imply 
$$Q_{j}(2i+1)=(-1)^{j}Q_{j}(2i)
\quad
\mbox{for $j=s-1$ and $s$}.$$

Equations (\ref{equ:ee}) and (\ref{equ:eo})
imply 
$$\sum_{j=0}^{[\frac{s}{2}]}Q_{2j}(2i)=
\sum_{j=0}^{[\frac{s-1}{2}]}Q_{2j+1}(2i)
=\delta_{i,0}\frac{|X|}{2}.$$
\begin{eqnarray}
&&Q_{2[\frac{s}{2}]}(2i)=\delta_{i,0}\frac{|X|}{2}
-\sum_{j=0}^{[\frac{s}{2}]-1}Q_{2j}(2i), \label{equ:4-12}\\
&&Q_{2[\frac{s-1}{2}]+1}(2i)=\delta_{i,0}\frac{|X|}{2}
-\sum_{j=0}^{[\frac{s-1}{2}]-1}Q_{2j+1}(2i)
\label{equ:4-13}
\end{eqnarray}
Therefore, if $s$ is even, then $\alpha_s=0$.
Hence $Q_{2j+1}(s)=\tilde Q_{2j+1}(\alpha_s)
=\tilde Q_{2j+1}(0)=0$ for any $j=0,\ldots,
\frac{s}{2}-2$. Hence (\ref{equ:4-13}) implies
\begin{equation}
Q_{s-1}(s)=-\sum_{j=0}^{\frac{s}{2}-2}Q_{2j+1}(s)=0.
\end{equation}
This completes the proof.\\

\noindent
{\bf Proof of Theorem 1.2 (3)}\\
If $i\leq s-2$, then $m_i=Q_i(0)
=\tilde Q_i(1)={n+i-1\choose i}
-{n+i-3\choose i-2}$.
If $s$ is odd, then (\ref{equ:4-12}) and (\ref{equ:4-13}) imply

\begin{eqnarray}
m_{s-1}&=&Q_{s-1}(0)=\frac{|X|}{2}
-\sum_{j=0}^{\frac{s-3}{2}}m_{2j}\nonumber\\
&=&\frac{|X|}{2}-\sum_{j=0}^{\frac{s-3}{2}}
\left({n+2j-1\choose 2j}-{n+2j-3\choose 2j-2}\right)
\nonumber\\
&=&\frac{|X|}{2}-{n+s-4\choose s-3}\\
m_s&=&Q_s(0)=\frac{|X|}{2}
-\sum_{j=0}^{\frac{s-3}{2}}m_{2j+1}\nonumber\\
&=&\frac{|X|}{2}
-\sum_{j=0}^{\frac{s-3}{2}}\left({n+2j\choose 2j+1}
-{n+2j-2\choose 2j-1}\right)
\nonumber\\
&=&\frac{|X|}{2}-{n+s-3\choose s-2}
\end{eqnarray}
If $s$ is even, then (\ref{equ:4-12}) and (\ref{equ:4-13}) imply
\begin{eqnarray}
m_{s-1}&=&Q_{s-1}(0)=\frac{|X|}{2}
-\sum_{j=0}^{\frac{s-4}{2}}m_{2j+1}\nonumber\\
&=&\frac{|X|}{2}-\sum_{j=0}^{\frac{s-4}{2}}
\left({n+2j\choose 2j+1}-{n+2j-2\choose 2j-1}\right)
\nonumber\\
&=&\frac{|X|}{2}-{n+s-4\choose s-3}\\
m_s&=&Q_s(0)=\frac{|X|}{2}
-\sum_{j=0}^{\frac{s-2}{2}}m_{2j}\nonumber\\
&=&\frac{|X|}{2}
-\sum_{j=0}^{\frac{s-2}{2}}\left({n+2j-1\choose 2j}
-{n+2j-3\choose 2j-2}\right)\nonumber\\
&=&\frac{|X|}{2}-{n+s-3\choose s-2}
\end{eqnarray}
This completes the proof.\\

\noindent
{\bf Proof of Theorem 1.2 (4)}\\ 
We compute $(|X|E_i\circ |X|E_j)\circ D_l$
in two ways. 
\begin{equation}(|X|E_i\circ |X|E_j)\circ D_l
=(\sum_{\mu=0}^sQ_i(\mu)D_\mu)\circ(\sum_{k=0}^sQ_j(k)D_k)\circ D_l=Q_i(l)Q_j(l)D_l
\label{equ:4-19}
\end{equation}
On the other hand we have
\begin{eqnarray}
&&(|X|E_i\circ |X|E_j)\circ D_l
=(\sum_{\mu=0}^sq_{i,j}^\mu|X|E_\mu)
\circ D_l
=(\sum_{\mu=0}^sq_{i,j}^\mu\sum_{k=0}^s
Q_\mu(k)D_k)\circ D_l\nonumber\\
&&=\sum_{\mu=0}^sq_{i,j}^\mu
Q_\mu(l) D_l.
\label{equ:4-20}
\end{eqnarray}
(\ref{equ:4-19}) and (\ref{equ:4-20}) imply
$$Q_i(l)Q_j(l)=\sum_{\mu=0}^s
q_{i,j}^\mu Q_\mu(l)$$
for any $0\leq l\leq s$. In particular 
we have
\begin{equation}Q_i(0)Q_j(0)=\sum_{\mu=0}^s
q_{i,j}^\mu Q_\mu(0).
\label{equ:3-1}
\end{equation}
Since Theorem \ref{theo:main2} (2) implies $Q_j(1)=(-1)^jQ_j(0)$,
we obtain
\begin{equation}Q_i(0)Q_j(0)=(-1)^{i+j}Q_i(1)Q_j(1)
=(-1)^{i+j}\sum_{\mu=0}^s
q_{i,j}^\mu Q_\mu(1)
=\sum_{\mu=0}^s(-1)^{i+j+\mu}
q_{i,j}^\mu Q_\mu(0).
\label{equ:3-2}
\end{equation}
Hence (\ref{equ:3-1}) and (\ref{equ:3-2})
imply
$$\sum_{\mu=0}^s( (-1)^{i+j+\mu}-1)
q_{i,j}^\mu Q_\mu(0)=0.$$
Since $Q_\mu(0)=m_\mu>0$ and
, $q_{i,j}^\mu\geq 0$ we must have 
$q_{i,j}^\mu=0$ for any $0\leq i,j,\mu\leq s$
satisfying $i+j+\mu$ is odd. 
This completes the proof.\\

\noindent
{\bf Proof of Theorem 1.2 (5)}\\
Now we are ready to proove that 
the association scheme given in Theorem 1.1
is Q-polynomial.
For $0\leq i\leq s-3$, we obtain
\begin{eqnarray}
&&|X|E_1\circ |X|E_i
=\sum_{l=0}^s\tilde Q_1(\alpha_l)D_l
\circ \sum_{\mu=0}^s\tilde Q_i(\alpha_\mu)D_\mu
=\sum_{l=0}^s\tilde Q_1(\alpha_l)
\tilde Q_i(\alpha_l)D_l\nonumber\\
&&=\sum_{l=0}^s\sum_{k=i-1, i+1}
q_k(1,i)
\tilde Q_k(\alpha_l)D_l
=\sum_{k=i-1,i+1}
q_k(1,i)E_k.
\nonumber
\end{eqnarray}
This implies
$q_{1,i}^{i-1}=q_{i-1}(1,i)=\frac{n(n+i-3)}{n+2i-4}$
for any $1\leq i\leq s-3$,
$q_{1,i}^{i+1}=q_{i+1}(1,i)=\frac{n(i+1)}{n+2i}$ 
for any $0\leq i\leq s-3$ and 
$q_{1,i}^j=0$ for any $i,\ j$ satisfying 
$0\leq i\leq s-3$ and
$0\leq j\neq i-1,i+1$.
Since Krein numbers satisfy 
$m_iq_{1,j}^i=m_jq_{1,i}^j$, and
$m_i>0$, we obtain $q_{1,s-1}^j=q_{1,s}^j=0$
for $0\leq j\leq s-3$
and $q_{1,s-2}^j=0$ for $0\leq j\leq s-4$.
Also (1) implies $q_{1,s-2}^s=q_{1,s}^{s-2}=0$.
Thus $B_1^*$ must be a tri-diagonal matrix
and the association scheme is Q-polynomial.

Next, we determin the nonzero entries of $B_1^*$.
Note that $F_{s-1}\in \langle E_{s-1},E_s\rangle$.
\begin{eqnarray}
&&|X|E_1\circ |X|E_{s-2}
=\sum_{l=0}^s\tilde Q_1(\alpha_l)D_l
\circ \sum_{\mu=0}^s\tilde Q_{s-2}(\alpha_\mu)D_\mu
=\sum_{l=0}^s\tilde Q_1(\alpha_l)
\tilde Q_{s-2}(\alpha_l)D_l\nonumber\\
&&=\sum_{l=0}^s\sum_{k=s-3, s-1}
q_k(1,s-2)
\tilde Q_k(\alpha_l)D_l
=|X|
q_{s-3}(1,s-2)E_{s-3}+|X|F_{s-1}.
\nonumber
\end{eqnarray}
Hence 
$q_{1,s-2}^{s-3}=q_{s-3}(1,s-2)=\frac{n(n+s-5)}
{n+2s-8}$.
The formula $\sum_{i=0}^sq_{1,i}^k=m_1(=n)$ for
any $0\leq k\leq s$ implies
$q_{1,s-1}^s=n$ and 
$q_{1,s-1}^{s-2}=n-q_{1,s-3}^{s-2}
=n-q_{s-2}(1,s-3)=\frac{n(n+s-4)}{n+2s-6}$.
Then we obtain
$$q_{1,s-2}^{s-1}=
\frac{m_{s-2}}{m_{s-1}}q_{1,s-1}^{s-2}
= \frac{2n(n-1)(n+s-4)!}{(s-2)!(n-1)!|X|-2(s-2)(n+s-4)!}. $$
Then 
$$q_{1,s}^{s-1}=n-q_{1,s-2}^{s-1}
= \frac{(s-2)!n!|X|-2n(n+s-3)!}
{(s-2)!(n-1)!|X|-2(s-2)(n+s-4)!}.$$
This completes the proof.
\hfill\qed\\

\section{Concluding Remarks}

We consider the case of $t\geq 2s-3$ and $s=4$, the smallest 
nontrivial case, since the case $s=3$ is a kind of studied as 
the case of equiangular lines in Euclidean spaces, (cf. \cite{L-S,
D-G-S}.) 
The details of this case of $t\geq 5$ and 
$s=4$ was already mentioned in our 
previous paper \cite{B-B-2}, so we just summarize the results. 
Note that the connection of this case with tight Euclidean $7$-designs is also mentioned in \cite{B-B-2}.\\

The first and the second eigen matrices $P$, $Q$ 
and the dual intersection 
matrix $B_1^*$, as well as all the intersection matrices 
$B_i(i=0,1,2,3,4)$ of the corresponding association scheme 
are described in Appendix 1. (Note that they are all 
described by two parameters $n$ and $N=|X|/2.$ )\\

We have shown that if $n\geq 3,$ then all the entries of 
$P$ and $Q$ are rational numbers. This gives a strong 
restriction on the possible pairs of $n$ and $N$ for which 
such association schemes may exist. It is also shown that 
for such an association scheme, if $A(X)=\{-1, -\alpha, 0, 
\alpha\}$ then $\alpha =\frac{1}{m}$ for a positive integer $m$. 
So, it is easy to see that if $m$ is given, then there are only 
finitely many possible pairs of $n$ and $N$. 
In \cite{B-B-2} 
we have listed the possible pairs $(n, N)$ for 
$m\leq 5.$ 
We 
have succeeded in 
classifying those with $m=2,3$ (see \cite{B-B-2}) and 
found some 
examples with $m=4$ and $N=144$. We recognized that 
this association scheme comes from the real MUB(mutually 
unbiased bases) in $\mathbb R^{16}$, and subsequently we noticed 
that such examples are obtained for any integer 
$m=\frac{1}{2^{2r}}$. 
Namely, we have a family of Q-polynomial association 
schemes $X$ with $|X|=2^{4r}+2^{2r+1}$ 
which gives the antipodal spherical design
in $S^{n-1}$, $n=2^{2r}$, of degree $4$ and 
strength $5$, with
 $A(X)=\{-1, -\frac{1}{2^{2r}}, 0,
\frac{1}{2^{2r}}\}$. 
These Q-polynomial association schemes are "not'' $P$-polynomial. 
(The explicit parameters are described in Appendix 2.) This family 
of association schemes was 
originally missing in the list of such association schemes in the home page of W. L. Martin 
(see \cite{M}). \\
We also mention that explicit descriptions of possible parameters with $s=5$ and $t=7$ are mentioned in 
\cite{B-B-2}. 
Anyway, it would be interesting to try to classify antipodal 
spherical $t$-designs with $t\geq 2s-3,$ in particular to show 
the nonexistence for large $t$ (and hence for large $s$.) \\

\noindent
{\bf Appendix 1.}\\
Eigenmatrices and the intersection matrices of the
Q-polynomial association schemes attached to the
spherical designs $X$ on $S^{n-1}$ of degree $4$ and strength $5$. We note that $|X|=2N$, $n\geq 3$,
$\frac{n(n+1)}{2}<N\leq \frac{n(n+1)(n+2)}{6}$ .\\

{\small \[B_0=
 \left[ \begin {array}{ccccc} 1&0&0&0&0\\\noalign{\medskip}0&1&0&0&0\\\noalign{\medskip}0&0&1&0&0\\\noalign{\medskip}0&0&0&1&0\\\noalign{\medskip}0&0&0&0&1\end {array} \right] 
 ,\qquad
 B_1=\left[ \begin {array}{ccccc} 0&1&0&0&0\\\noalign{\medskip}1&0&0&0&0\\\noalign{\medskip}0&0&0&1&0\\\noalign{\medskip}0&0&1&0&0\\\noalign{\medskip}0&0&0&0&1\end {array} \right] ,
 \]}
 \noindent
{\small
$B_2=$
  \[\left[ \begin {array}{ccc} 
 0&0&1\\
 \noalign{\medskip}
 0&0&0\\
 \noalign{\medskip}
 {\frac { \left( N-n \right) ^{2} \left( n+2 \right) }{n \left( 3\,N-{n}^{2}-2\,n \right) }}&
 0&
 {\frac { \left( N-n \right)  \left( {N}^{2}n+8\,{N}^{2}-9\,{n}^{2}N-18\,nN+2\,{n}^{4}+8\,{n}^{3}+8\,{n}^{2} \right) }{2n \left( 3\,N-{n}^{2}-2\,n \right) ^{2}}}
 - \frac{ N-2n}{2n}\sqrt {{\frac { \left( n+2 \right)  \left( N-n \right) }{3\,N-{n}^{2}-2\,n}}} \\
 
 \noalign{\medskip}0&
 {\frac { \left( N-n \right) ^{2} \left( n+2 \right) }{n \left( 3\,N-{n}^{2}-2\,n \right) }}&
 \frac{ N-2\,n } {2n}\sqrt {{\frac { \left( n+2 \right)  \left( N-n \right) }{3\,N-{n}^{2}-2\,n}}} 
 
 +{\frac { \left( N-n \right)  \left( {N}^{2}n+8\,{N}^{2}-9\,{n}^{2}N-18\,nN+2\,{n}^{4}+8\,{n}^{3}+8\,{n}^{2} \right) }{2n \left( 3\,N-{n}^{2}-2\,n \right) ^{2}}}\\
 
 \noalign{\medskip}0&
 0&
 {\frac {{N}^{2} \left( n-1 \right)  \left( 2\,N-{n}^{2}-n \right) }{n \left( 3\,N-{n}^{2}-2\,n \right) ^{2}}}\end {array} \right.
\]

 \[
  \left. \begin {array}{cc} 
 0&0\\
 \noalign{\medskip}
 1&0\\
 \noalign{\medskip}
  
\frac{ N-2\,n } {2n} \sqrt {{\frac { \left( n+2 \right)  \left( N-n \right) }{3\,N-{n}^{2}-2\,n}}} +{\frac { \left( N-n \right)  \left( {N}^{2}n+8\,{N}^{2}-9\,{n}^{2}N-18\,nN+2\,{n}^{4}+8\,{n}^{3}+8\,{n}^{2} \right) }{2n \left( 3\,N-{n}^{2}-2\,n \right) ^{2}}}&
 
 {\frac {N \left( -N+n \right) ^{2} \left( n+2 \right) }{2n \left( -3\,N+{n}^{2}+2\,n \right) ^{2}}}\\
 
 \noalign{\medskip}
 
 {\frac { \left( N-n \right)  \left( {N}^{2}n+8\,{N}^{2}-9\,{n}^{2}N-18\,nN+2\,{n}^{4}+8\,{n}^{3}+8\,{n}^{2} \right) }{2n \left( 3\,N-{n}^{2}-2\,n \right) ^{2}}}
 
 -
\frac{ N-2\,n }{2n} \sqrt {{\frac { \left( n+2 \right)  \left( N-n \right) }{3\,N-{n}^{2}-2\,n}}} &

 {\frac {N \left( N-n \right) ^{2} \left( n+2 \right) }{2n \left( 3\,N-{n}^{2}-2\,n \right) ^{2}}}
 \\
 
 \noalign{\medskip}

 {\frac {{N}^{2} \left( n-1 \right)  \left( 2\,N-{n}^{2}-n \right) }{n \left( 3\,N-{n}^{2}-2\,n \right) ^{2}}}&
 
 {\frac { \left( N-n \right) ^{2} \left( n+2 \right)  \left( 2\,N-{n}^{2}-2\,n \right) }{n \left( 3\,N-{n}^{2}-2\,n \right) ^{2}}}\end {array} \right]
 \]}

\newpage
\noindent
$B_3=$\\
{\small \[
\left[ \begin {array}{ccc} 
0&0&0\\
\noalign{\medskip}
0&0&1\\
\noalign{\medskip}0&
{\frac { \left( N-n \right) ^{2} \left( n+2 \right) }{n \left( 3\,N-{n}^{2}-2\,n \right) }}&
\frac{ N-2\,n }{2n}
\sqrt {{\frac { \left( n+2 \right)  \left( N-n \right) }{3\,N-{n}^{2}-2\,n}}} 

+{\frac { \left( N-n \right)  \left( {N}^{2}n+8\,{N}^{2}-9\,{n}^{2}N-18\,nN+2\,{n}^{4}+8\,{n}^{3}+8\,{n}^{2} \right) }{2n \left( 3\,N-{n}^{2}-2\,n \right) ^{2}}}\\

\noalign{\medskip}
{\frac { \left( N-n \right) ^{2} \left( n+2 \right) }{n \left( 3\,N-{n}^{2}-2\,n \right) }}&
0&
{\frac { \left( N-n \right)  \left( {N}^{2}n+8\,{N}^{2}-9\,{n}^{2}N-18\,nN+2\,{n}^{4}+8\,{n}^{3}+8\,{n}^{2} \right) }{2n \left( 3\,N-{n}^{2}-2\,n \right) ^{2}}}
-
\frac{ N-2\,n }{2n}\sqrt {{\frac { \left( n+2 \right)  \left( N-n \right) }{3\,N-{n}^{2}-2\,n}}} \\

\noalign{\medskip}
0&
0&
{\frac {{N}^{2} \left( n-1 \right)  \left( 2\,N-{n}^{2}-n \right) }{n \left( 3\,N-{n}^{2}-2\,n \right) ^{2}}}
\end {array} \right.
\]

\[\left.
\begin {array}{cc} 1&0\\
\noalign{\medskip}
0&0\\
\noalign{\medskip}

{\frac { \left( N-n \right)  \left( {N}^{2}n+8\,{N}^{2}-9\,{n}^{2}N-18\,nN+2\,{n}^{4}+8\,{n}^{3}+8\,{n}^{2} \right) }{2n \left( 3\,N-{n}^{2}-2\,n \right) ^{2}}}

-\frac{ N-2\,n }{2n}
\sqrt {{\frac { \left( n+2 \right)  \left( N-n \right) }{3\,N-{n}^{2}-2\,n}}} &
{\frac {N \left( N-n \right) ^{2} \left( n+2 \right) }{2n \left( 3\,N-{n}^{2}-2\,n \right) ^{2}}}\\

\noalign{\medskip}
\frac{N-2\,n}{2n}
\sqrt {{\frac { \left( n+2 \right)  \left( N-n \right) }{3\,N-{n}^{2}-2\,n}}} 

+{\frac { \left( N-n \right)  \left( {N}^{2}n+8\,{N}^{2}-9\,{n}^{2}N-18\,nN+2\,{n}^{4}+8\,{n}^{3}+8\,{n}^{2} \right) }{2n \left( 3\,N-{n}^{2}-2\,n \right) ^{2}}}&
{\frac {N \left( N-n \right) ^{2} \left( n+2 \right) }{2n \left( 3\,N-{n}^{2}-2\,n \right) ^{2}}}\\

\noalign{\medskip}

{\frac {{N}^{2} \left( n-1 \right)  \left( 2\,N-{n}^{2}-n \right) }{n \left( 3\,N-{n}^{2}-2\,n \right) ^{2}}}&
{\frac { \left( N-n \right) ^{2} \left( n+2 \right)  \left( 2\,N-{n}^{2}-2\,n \right) }{n \left( 3\,N-{n}^{2}-2\,n \right) ^{2}}}\end {array} \right]
\]}

{\small
\[B_4=
\left[ \begin {array}{ccc} 
0&0&0\\
\noalign{\medskip}
0&0&0\\
\noalign{\medskip}
0&0&
{\frac {{N}^{2} \left( n-1 \right)  \left( 2\,N-{n}^{2}-n \right) }{n \left( 3\,N-{n}^{2}-2\,n \right) ^{2}}}
\\
\noalign{\medskip}
0&
0&
-{\frac {{N}^{2} \left( n-1 \right)  \left( -2\,N+{n}^{2}+n \right) }{n \left( -3\,N+{n}^{2}+2\,n \right) ^{2}}}\\

\noalign{\medskip}
{\frac {2\,N \left( n-1 \right)  
\left( 2\,N-{n}^{2}-n \right) }
{n \left( 3\,N-{n}^{2}-2\,n \right) }}&
{\frac {2\,N \left( n-1 \right)  \left( -2\,N+{n}^{2}+n \right) }{n \left( -3\,N+{n}^{2}+2\,n \right) }}&
{\frac {2\,N \left( n-1 \right) 
 \left( 2\,N-{n}^{2}-n \right)  \left( 2\,N-{n}^{2}-2\,n \right) }{n \left( 3\,N-{n}^{2}-2\,n \right) ^{2}}}
\end {array} \right.
\]

\[
\left. \begin {array}{cc} 
0&1\\

\noalign{\medskip}
0&1\\
\noalign{\medskip}
{\frac {{N}^{2} \left( n-1 \right) 
 \left( 2\,N-{n}^{2}-n \right) }
 {n \left( 3\,N-{n}^{2}-2\,n \right) ^{2}}}&
{\frac { \left( n+2 \right)  \left( N-n \right) ^{2}
 \left(2\,N-{n}^{2}-2\,n \right) }
 {n \left( 3\,N-{n}^{2}-2\,n \right) ^{2}}}
\\
\noalign{\medskip}

{\frac {{N}^{2} \left( n-1 \right) 
 \left(2\,N-{n}^{2}-n \right) }
 {n \left( 3\,N-{n}^{2}-2\,n \right) ^{2}}}&
{\frac { \left( n+2 \right) 
 \left( N-n \right) ^{2} 
 \left(2\,N-{n}^{2}-2\,n \right) }
 {n \left( 3\,N-{n}^{2}-2\,n \right) ^{2}}}
\\
\noalign{\medskip}

{\frac {2\,N \left( n-1 \right) 
 \left( 2\,N-{n}^{2}-n \right) 
  \left( 2\,N-{n}^{2}-2\,n \right) }
  {n \left( 3\,N-{n}^{2}-2\,n \right) ^{2}}}&
{\frac {2\, \left( 4\,{N}^{2}n-4\,{n}^{3}N+6\,{n}^{2}N+{n}^{5}-5\,{n}^{3}-10\,{N}^{2}+10\,nN-2\,{n}^{2} \right) N}
{n \left( 3\,N-{n}^{2}-2\,n \right) ^{2}}}\end {array} \right]
\]}
\[P=
\left[ \begin {array}{ccccc} 
1&1&
{\frac { \left( N-n \right) ^{2} \left( n+2 \right) }{n \left( 3\,N-{n}^{2}-2\,n \right) }}&
{\frac { \left( N-n \right) ^{2} \left( n+2 \right) }{n \left( 3\,N-{n}^{2}-2\,n \right) }}&
{\frac {2\, \left( n-1 \right) N \left( 2\,N-{n}^{2}-n \right) }{n \left( 3\,N-{n}^{2}-2\,n \right) }}\\
\noalign{\medskip}
1&-1&
- \frac{ N-n }{n}\sqrt {{\frac { \left( n+2 \right)  \left( N-n \right) }
{3\,N-{n}^{2}-2\,n}}}& 
\frac{N-n }{n}\sqrt {{\frac { \left( n+2 \right)  \left( N-n \right) }{3\,N-{n}^{2}-2\,n}}}&
0\\\noalign{\medskip}1&
1&
{\frac { \left( N-n \right)  \left( 2\,N-{n}^{2}-2\,n \right) }{n \left( 3\,N-{n}^{2}-2\,n \right) }}&
{\frac { \left( N-n \right)  \left( 2\,N-{n}^{2}-2\,n \right) }{n \left( 3\,N-{n}^{2}-2\,n \right) }}&
-{\frac {2\,N \left( 2\,N-{n}^{2}-n \right) }{n \left( 3\,N-{n}^{2}-2\,n \right) }}\\
\noalign{\medskip}1&-1&\sqrt {{\frac { \left( n+2 \right)  \left( N-n \right) }{3\,N-{n}^{2}-2\,n}}}&-\sqrt {{\frac { \left( n+2 \right)  \left( N-n \right) }{3\,N-{n}^{2}-2\,n}}}&0\\
\noalign{\medskip}1&1&-{\frac { \left( n+2 \right)  \left( N-n \right) }{3\,N-{n}^{2}-2\,n}}&-{\frac { \left( n+2 \right)  \left( N-n \right) }{3\,N-{n}^{2}-2\,n}}&
{\frac {2\, \left( n-1 \right) N}{3\,N-{n}^{2}-2\,n}}\end {array} \right] 
\]

\[Q=
\left[ \begin {array}{ccccc} 1&n&
\frac{ 
\left( n+2 \right)  \left( n-1 \right)}{2} 
&N-n&
\frac{2\,N-{n}^{2}-n}{2}
\\
\noalign{\medskip}1&-n&\frac{\left( n+2 \right)  
\left( n-1 \right)}{2} &
-(N-n)&
\frac{2\,N-{n}^{2}-n}{2}\\
\noalign{\medskip}
1&
-n{{\sqrt {{\frac {3\,N-{n}^{2}-2\,n}{ \left( n+2 \right)  
\left(N-n \right) }}}}}&
{\frac {\left( n-1 \right) 
 \left( 2\,N-{n}^{2}-2\,n \right) }
 {2\, \left(N-n\right)}}&
n{{\sqrt {{\frac {3\,N-{n}^{2}-2\,n}{ \left( n+2 \right)  
\left( N-n \right) }}}}}
&
-{\frac { \left( 2\,N-{n}^{2}-n \right) n}{2\,\left(N-n\right)}}\\
\noalign{\medskip}
1&n{{\sqrt {{\frac {3\,N-{n}^{2}-2\,n}{ \left( n+2 \right) 
 \left( N-n \right) }}}}}&
{\frac { \left( n-1 \right)  \left( 2\,N-{n}^{2}-2\,n \right) }
{2\,\left(N-n\right)}}&
-n{{\sqrt {{\frac {3\,N-{n}^{2}-2\,n}
{ \left( n+2 \right)  \left( N-n \right) }}}}}&
-{\frac { \left( 2\,N-{n}^{2}-n \right) n}
{2\,\left(N-n\right)}}\\\noalign{\medskip}1&0&
-\frac{n+2}{2}&
0&
\frac{n}{2}\end {array} \right] 
\]

\[
B_1^*= 
\left[ \begin {array}{ccccc} 0&1&0&0&0\\
\noalign{\medskip}n&0&{\frac {2\,n}{n+2}}&0&0\\
\noalign{\medskip}0&n-1&0&
{\frac { \left( n-1 \right) {n}^{2}}{2\,\left(N-n\right)}}&
0\\
\noalign{\medskip}0&0&{\frac {{n}^{2}}{n+2}}&0&n\\\noalign{\medskip}0&0&0&{\frac { \left( 2\,N-{n}^{2}-n \right) n}
{2\,\left(N-n\right)}}&0\end {array} \right] \]

\noindent
{\bf Appendix 2.}\\
Association schemes attached to the spherical 
design $X$ in $S^{n-1}$ of degree $4$ and strength $5$ 
with $n=2^{2r}, r\geq 1$ and $|X|=2N=2^{4r}+2^{2r+1}$\\

\begin{eqnarray}
B_0&=&\left[\begin{array}{ccccc}
1&0&0&0&0\\
0&1&0&0&0\\
0&0&1&0&0\\
0&0&0&1&0\\
0&0&0&0&1
\end{array}
\right],\qquad 
B_1=\left[\begin{array}{ccccc}
0&1&0&0&0\\
1&0&0&0&0\\
0&0&0&1&0\\
0&0&1&0&0\\
0&0&0&0&1
\end{array}
\right],\nonumber\\
\noalign{\medskip}
B_2&=&\left[\begin{array}{ccccc}
0&0&1&0&0\\
0&0&0&1&0\\
2^{4r-1}&0&2^{r-2}(2^{2r}-2)(2^r-1)
&2^{r-2}(2^{2r}-2)(2^r+1)&2^{4r-1}\\
0&2^{4r-1}&2^{r-2}(2^{2r}-2)(2^r+1)
&2^{r-2}(2^{2r}-2)(2^r-1)
&2^{4r-1}\\
0&0&2^{2r}-1&2^{2r}-1&0
\end{array}
\right]\nonumber
\end{eqnarray}
\begin{eqnarray}
B_3&=&\left[\begin{array}{ccccc}
0&0&0&1&0\\
0&0&1&0&0\\
0&2^{4r-1}
&2^{r-2}(2^{2r}-2)(2^r+1)
&2^{r-2}(2^{2r}-2)(2^r-1)
&2^{4r-1}\\
2^{4r-1}&0
&2^{r-2}(2^{2r}-2)(2^r-1)
&2^{r-2}(2^{2r}-2)(2^r+1)&2^{4r-1}\\
0&0&2^{2r}-1&2^{2r}-1&0
\end{array}
\right],\nonumber\\
\noalign{\medskip}
B_4&=&\left[\begin{array}{ccccc}
0&0&0&0&1\\
0&0&0&0&1\\
0&0&2^{2r}-1&2^{2r}-1&0\\
0&0&2^{2r}-1&2^{2r}-1&0\\
2(2^{2r}-1)&2(2^{2r}-1)&0&0&2(2^{2r}-2)
\end{array}
\right],\nonumber\\
\noalign{\medskip}
P&=&\left[\begin{array}{ccccc}
1&1&2^{4r-1}&2^{4r-1}&2(2^{2r}-1)\\
1&-1&-2^{3r-1}&2^{3r-1}&0\\
1&1&0&0&-2\\
1&-1&2^r&-2^r&0\\
1&1&-2^{2r}&-2^{2r}&2(2^{2r}-1)
\end{array}\right],\nonumber\\
\noalign{\medskip}
Q&=&\left[\begin{array}{ccccc}
1&2^{2r}&(2^{2r}-1)(2^{2r-1}+1)
&2^{4r-1}&2^{2r-1}\\
1&-2^{2r}&(2^{2r}-1)(2^{2r-1}+1)&-2^{4r-1}&2^{2r-1}\\
1&-2^r&0&2^r&-1\\
1&2^r&0&-2^r&-1\\
1&0&-2^{2r-1}-1&0&2^{2r-1}
\end{array}\right],\nonumber\\
\noalign{\medskip}
B_1^*&=&\left[\begin{array}{ccccc}
0&1&0&0&0\\
2^{2r}&0&\frac{2^{2r}}{2^{2r-1}+1}&0&0\\
0&2^{2r}-1&0&2^{2r}-1&0\\
0&0&\frac{2^{4r}}{2^{2r}+2}&0&2^{2r}\\
0&0&0&1&0
\end{array}\right]\nonumber
\end{eqnarray}

\newpage
\noindent
2000 Mathematics
Subject Classification.

Primary: 05E99, Secondary:
05B99, 51M99, 62K99.\\

\noindent
Key words and Phrases.

spherical design, 
tight design, 
$s$-distance set, association scheme

\end{document}